\input amssym.def 
\input amssym
\magnification=1200
\parindent0pt
\hsize=16 true cm
\baselineskip=13  pt plus .2pt
$ $

\centerline {\bf  On geodesic embeddings of hyperbolic surfaces}
\centerline {\bf into hyperbolic 3-manifolds}

\bigskip

\centerline {Bruno P. Zimmermann}

\bigskip \bigskip

{\bf Abstract.}  We consider the problem of when a closed orientable hyperbolic
surface admits a totally geodesic embedding into a closed orientable hyperbolic
3-manifold; given  a finite isometric group action
on the surface, we consider in particular equivariant versions of such an embedding.
We prove that an equivariant embedding exists for all finite {\it irreducible} group
actions on surfaces; such surfaces are known also as
{\it quasiplatonic surfaces}; in particular, all quasiplatonic surfaces embed
geodesically. In the last section, we discuss some cases of more general finite
group actions on surfaces.

\medskip

e-mail:   zimmer@units.it

\bigskip \bigskip

{\bf 1. Introduction}

\medskip

In the present paper, all surfaces and 3-manifolds will be compact, connected and
orientable, all group-actions orientation-preserving. Let $\cal S$ be a closed 
hyperbolic surface and  $G$ a finite group of isometries of $\cal S$; we say that the
pair $({\cal S},G)$ {\it embeds geodesically} 
if there exists a totally geodesic embedding of 
$\cal S$ into a closed hyperbolic 3-manifold ${\cal M}$ such that the action of $G$ on
${\cal S}$ extends to an isometric action of $G$ on ${\cal M}$; if $G$ is the trivial
group, we just say that the hyperbolic surface ${\cal S}$ embeds geodesically.
Moreover, the pair
$({\cal S},G)$ {\it bounds geometrically} if ${\cal S}$ is the unique, totally geodesic 
boundary component of a compact hyperbolic 3-manifold ${\cal M}$ and the action of
$G$ on ${\cal S}$ extends to an isometric action of $G$ on ${\cal M}$.

\medskip

It is shown in [GZ] (see also the survey [Z1])
that various infinite series of Hurwitz-actions of maximal possible order $84(g-1)$ on
a surface of genus $g$ bound geometrically; in particular, the underlying
Hurwitz surfaces bound geometrically. Finite group actions on surfaces of small
genus are considered in  [WZ] and [Z2]; for example, the smallest Hurwitz action of
order 168 on Klein's quartic of genus 3 does not bound any compact 3-manifold with a
unique boundary component but it remains open whether Klein's quartic bounds
geometrically; however, as a consequence of  [Z3, Proposition 1], every Hurwitz action
embeds geodesically.

\medskip

Lifting a finite group $G$ of isometries of a hyperbolic surface ${\cal S}$ to its
universal covering $\Bbb H^2$ we get a Fuchsian group $F$ and an exact sequence 
$1 \to \pi_1({\cal S})  \to F \to G \to 1$ (where $\pi_1({\cal S})$ is interpreted as
the group of covering transformations). We are interested here in the cases where $F$
is a  triangle group $(p,q,r)$ (i.e., the orientation-preserving subgroup of index 2 of
the group generated by the reflections in the sides of a hyperbolic triangle with
angles 
$\pi/p, \pi/q$ and $\pi/r$), with a presentation 
$$<x, y  \mid  x^p = y^q = (xy)^r = 1>.$$ 
These are exactly the {\it irreducible actions} of finite groups on surfaces (cf. 
[Z1]), and we say that the action of $G$ is of triangle type 
$(p,q,r)$. Since any two triangle groups of the same type are conjugate by an isometry, 
an irreducible action of a finite group $G$ determines the hyperbolic structure of the
underlying surface
${\cal S}$ (this is true only  for Fuchsian triangle groups since the Teichm\"uller
or deformation space of any other Fuchsian group has positive dimension).

\medskip

We call a hyperbolic surface $\cal S$ {\it semiplatonic} if its universal
covering group is a subgroup of finite index in a triangle group $(p,q,r)$;    
if the universal covering is a {\it normal} subgroup of a triangle group $(p,q,r)$ then
the surface
$\cal S$ admits an irreducible action of the factor group
$G = (p,q,r)/\pi_1(\cal S)$, and such surfaces are known also as {\it quasiplatonic
surfaces}  (cf. [SW]).

\bigskip

{\bf Theorem 1.}   {\sl   Every finite irreducible group action $({\cal S},G)$  of a
finite group $G$ on a closed  hyperbolic surface ${\cal S}$  embeds geodesically into
a closed hyperbolic 3-manifold; in particular, every quasiplatonic surface 
${\cal S}$ embeds geodesically.}

\bigskip

{\bf Corollary 1.}   {\sl Every Hurwitz action of maximal possible order $84(g-1)$ on a
closed hyperbolic surface of genus $g$ embeds geodesically into a closed hyperbolic
3-manifold; in particular, all Hurwitz surfaces embed geodesically. }

\bigskip

However, as noted before, not all Hurwitz actions bound geometrically (not even
topologically, see [Z1]).

\bigskip

{\bf Corollary 2.}   {\sl Every semiplatonic hyperbolic surface embeds
geodesically into a closed hyperbolic 3-manifold.}

\bigskip

In general, finite {\it reducible actions} on surfaces do not embed geodesically since
for each topological type of such an action there is an
uncountable continuum of geometric realizations $({\cal S},G)$  but only countably
many can embed geometrically.  However we state the following:

\medskip

{\bf Conjecture.}  Let $G$ be a finite group of homeomorphisms of a
closed surface
$\cal S$ of genus $g > 1$. Then some geometric realization of $({\cal S},G)$ by a
hyperbolic surface $\cal S$ and a group of isometries $G$ of $\cal S$ embeds
geodesically into a closed hyperbolic 3-manifold.

\medskip

For irreducible actions, there is a unique geometric realization $({\cal S},G)$ which,
by Theorem 1, embeds geodesically. In the next theorem we consider the case of the
reducible actions which realize the largest possible orders after the irreducible
ones.  Let $G$ be a finite group of homeomorphisms of a closed surface of genus $g
> 1$. By choosing a hyperbolic structure of the quotient 2-orbifold  
${\cal S}/G$ and lifting to $\cal S$, one realizes $\cal S$ as a
hyperbolic surface on which $G$ acts by isometries. By lifing $G$ to
the universal covering $\Bbb H^2$ one obtains a Fuchsian group, and we consider the
case of a {\it quadrangle group} $(p,q,r,s)$ which has a presentation 
$$<x, y ,z, w \mid  x^p = y^q = z^r = w^s = xyzw = 1>;$$ 
we say that the action of $G$ is of {\it quadrangle type} $(p,q,r,s)$.

\bigskip

{\bf  Theorem 2.}   {\sl  Let $G$ be a finite group of quadrangle type of homeomorphisms
of a surface
$\cal S$ of genus $g > 1$. Then, for some geometric realization of 
$({\cal S},G)$ by a hyperbolic surface $\cal S$ and an isometric action of $G$, the
pair $({\cal S},G)$ embeds geodesically into a closed hyperbolic 3-manifold.}

\bigskip

In section 3 we present some results for other types of finite group
actions and Fuchsian groups.

\bigskip \bigskip

{\bf 2. Proof of  Theorem 1}

\medskip

i)  We give the proof first for the case of an irreducible action of
triangle type $(2,q,r)$, and then describe the changes necessary for the general case
$(p,q,r)$.

\medskip

Let $({\cal S},G)$ be a finite irreducible group action of triangle type
$(2,q,r)$.  Let ${\cal T}(2,q,r)$ be a truncated hyperbolic tetrahedron 
with dihedral angles $\pi/2, \pi/q$ and $\pi/r$ (or with labels 2, $q$ and $r$) at the
three pairs of opposite edges, truncated by orthogonal planes at the four virtual
vertices (so  ${\cal T}(2,q,r)$ has four hexagonal faces,  and four triangular
faces on the truncating planes).   Let $T(2,q,r)$
denote the tetrahedral group associated to ${\cal T}(2,q,r)$,
i.e. the orientation-preserving subgroup of index 2 in the Coxeter group generated by
the reflections in the four hexagonal faces of the tetrahedron. Let $a, b, x, y$
denote the rotations of orders $q, 2, 2, q$ around the edges with these labels
(suitably oriented). We represent the six labelled edges of the tetrahedron by a
square and its two
diagonals, avoiding the intersection point of the diagonals by an
over/under-crossing and removing small open neighborhoods of the four vertices, with
label $r$ at the two upper and lower edges of the square, label $q$ at the right and
the left edge, and label 2 at the two diagonals; then, using standard methods for
the computation of an orbifold fundamental group, one obtains the following
presentation of the tetrahedral group  $T(2,q,r)$ (see  [Z4] or [GZ] for similar
computations; see also [B]):

$$<a,b,x,y  \mid  b^q = a^2 = x^2 = y^q = 1, \; ab = xy, \; (ab)^r = (xy)^r =
(x^{-1}b)^r = 1 > $$  
(with the rotations $b$ and $y$ suitably oriented; also, $((x^{-1}a^{-1}x)y)^r = 1$).
The sugroup of
$T(2,q,r)$ generated by $x$ and  $y$ is a triangle group of type $(2,q,r)$, with the
presentation given in section 1 (with $p = 2$), and the surjection 
$\phi: (2,q,r) \to  G$ associated to the action of $G$ on ${\cal S}$ extends to a
surjection $\psi: T(2,q,r) \to G$ by setting  $\psi(x) = \phi(x)$, $\psi(y) = \phi(y)$, 
$\psi (a) = \phi(x)$ and $\psi(b) = \phi(y)$; note that $\psi((x^{-1}b)^r) =  
(\psi(x^{-1})\psi(b))^r = (\psi(x)\psi(b))^r = (\phi(x)\phi(y))^r = \phi(xy)^r = 1$, so
this is exactly the point where we use that $x$ has order 2.

\medskip

Let $K  \subset T(2,q,r)$ denote the kernel of $\psi$.  Let $\tilde \Bbb H^3$ denote
the hyperbolic 3-space truncated by all images under $T(2,q,r)$ of the four truncating
planes of the tetrahedron ${\cal T}(2,q,r)$. The quotient $M = \tilde
\Bbb H^3/K$  is a compact hyperbolic 3-manifold $\cal M$ with an isometric action of $G$
(the projection of $T(2,q,r)$), and  with totally geodesic boundary consisting of four
copies of the hyperbolic surface
${\cal S}$,  each with an action of type $(2,q,r)$ of $G$ (all isometrically
equivalent).  To obtain a closed hyperbolic 3-manifold with a geodesic embedding of
$({\cal S},G)$, one can either take the double of ${\cal M}$ and its $G$-action along
the boundary of $\cal M$, or alternatively glue in pairs the four boundary surfaces of
${\cal M}$ by
$G$-equivariant  isometries.

\medskip

This completes the proof of the Theorem 1 for the case of an action of type $(2,q,r)$.
The general case
$(p,q,r)$ is similar but we replace the hyperbolic 3-orbifold 
$\tilde \Bbb H^3/K$ by a  more complicated one and indicate the necessary
changes in the following  (see also the Remark in section 3 for another proof).

\medskip

The underlying topological space of the hyperbolic 3-orbifold ${\cal O} =
\tilde \Bbb H^3/K$ (compact with four totally geodesic boundary components) is the
3-sphere minus the interiors of four disjoint 3-balls, its singular set  consists of
the three pairs of opposite edges with labels $2, q$ and $r$ of the truncated
tetrahedron ${\cal T}(2,q,r)$. Representing the singular set by a square with diagonals
as before, we  associate label $p$ now to the two diagonals (instead of label 2 as 
before). This defines a compact hyperbolic 3-orbifold ${\cal O'}$ with four totally
geodesic boundary components, each a 2-sphere with three singular points of orders
$p$, $q$ and $r$. The two diagonals have one over/under crossing or half-twist which
we replace by $u$ half-twists  (for an integer $u$); this defines a 3-orbifold
${\cal O}(u)$ which is obtained by orbifold Dehn surgery on the orbifold ${\cal O'}$
(see  [Z4] or [GZ] for some figures).
By the orbifold version of Thurston's hyperbolic Dehn surgery theorem, for large enough
values of $u$ the orbifold ${\cal O}(u)$ is hyperbolic again, compact and with totally
geodesic boundary (in fact, we believe that $u = 2$ is sufficient, applying the
orbifold version of Thurston's hyperbolization theorem).

\medskip

As before, the orbifold fundamental group $\pi_1({\cal O}(u))$ is generated by
elements $a,b,x,y$ (see [Z4], [GZ]).  We suppose now that $u$ is even; then, as in the
first part of the proof, the surjection $\phi: (p,q,r)  \to G$ associated to the action
of $G$ on ${\cal S}$ extends to a surjection $\psi: \pi_1({\cal O}(u)) \to G$  (which 
does not work if $u$ is odd unless $p = 2$). Now the proof is completed exactly as in
the first case. 

\medskip

For an alternative proof in the second case, see the Remark in the next section.

\bigskip 
\vfill  \eject

{\it Proof of Corollary 2}.  For a closed hyperbolic surface ${\cal S}$, suppose that
its universal covering group $\pi_1({\cal S})$ is a subgroup of finite index in a
triangle group $(p,q,r)$. By taking the intersection of $\pi_1({\cal S})$ with all of
its conjugates in $(p,q,r)$, one obtains a  normal subgroup of finite index of
$(p,q,r)$ which is the fundamental group $\pi_1(\tilde {\cal S})$ of a finite
covering $\tilde {\cal S}$ of ${\cal S}$. Then $\tilde {\cal S}$ is a quasiplatonic
surface of type $(p,q,r)$, with an isometric action of $G = (p,q,r)/\pi_1(\tilde {\cal
S})$ associated to a surjection $\phi: (p,q,r) \to G$ with kernel 
$\pi_1(\tilde {\cal S})$.  Let $G_0$ denote the image
$\phi(\pi_1({\cal S}))$ of $\pi_1({\cal S}) \subset (p,q,r)$ in $G$; $G_0$ is the
covering group of the covering 
$\tilde {\cal S}$ of ${\cal S}$,  with $\tilde {\cal S}/G_0 =  {\cal S}$.

\medskip

By Theorem 1, $(\tilde {\cal S},G)$ embeds
geodesically into a closed hyperbolic 3-manifold $\cal M$ with an isometric
$G$-action, and ${\cal M}/G_0$ is a closed hyperbolic 3-manifold with a totally geodesic
embedding of the surface $\tilde {\cal S}/G_0 = {\cal S}$ (note that, by definition
of $\psi$ in the proof of Theorem 1, $\psi^{-1}(G_0)$ is torsionfree).

\bigskip \bigskip

{\bf 3. Proof of  Theorem 2}

\medskip

As before, the action of $G$ on $\cal S$ is defined by a surjection  
$\phi: (p,q,r,s) \to  G$ with kernel $\pi_1(\cal S)$, and we have the presentation 
$<x, y ,z, w \mid  x^p = y^q = z^r = w^s = xyzw = 1>$ of the quadrangle group  
$(p,q,r,s)$.

\medskip

We consider a hyperbolic 3-orbifold ${\cal O}$ with totally geodesic boundary whose
underlying topological space is $S^3$ minus the interiors of 6 small disjoint 3-balls
around the vertices of a standardly embedded octahedron, and whose singular set
$\Sigma$ is the 1-skeleton of the octahedron, truncated at its six vertices by the six
3-balls.  We assign labels $p,q,r$ and $s$ to the 12 edges of the
octahedron in the following way.

\medskip

We consider a planar projection of the 1-skeleton of the
octahedron, consisting of a square with its diagonals representing 8 edges of the
octahdron, a vertex outside of the square connected by four further
edges to the vertices of the square, and also a direction of
undercrossing for each of the 12 (disjoint) edges. It is easy to verify then that one
can assign labels $p,q,r$ and $s$ to the 12 edges such that at four of the vertices
of the projected 1-skeleton of the octahedron, the labels are in cyclic order $p,q,r$
and $s$ (following the chosen directions of the undercrossings), and at each of the
remaining two vertices the four edges have the same label, two times in the chosen and
two times in the opposite direction.

\medskip

By Andreev's theorem on hyperbolic polyhedra ([A], [T, chapter 13], [V, p.111]), the
truncated octahedron can be realized in $\Bbb H^3$ as a hyperbolic, orthogonally
truncated octahedron,  with dihedral angles
$\pi/p, \pi/q, \pi/r$ and $\pi/s$ according to the labels of the 12 edges (since there
are no essential Euclidean 2-suborbifolds).  Let 
$O$ denote the orientation-preserving subgroup of index 2 in the group generated by the
reflections in the faces of this hyperbolic octahedron. Then $O$ acts on 
$\bar \Bbb H^3$, obtained by truncating $\Bbb H^3$ by the truncating planes of the
octahedron and its images under $O$, and the quotient-orbifold  $\bar \Bbb H^3/O$
realizes the hyperbolic 3-orbifold $\cal O$, with orbifold fundamental group 
$\pi_1({\cal O}) \cong O$.

\medskip

A presentation of the orbifold fundamental group $\pi_1(\cal O)$ of $\cal O$ can be
obtained in a standard way from the planar projection of the singular set 
$\Sigma$ (similar as the Wirtinger presentation of the fundamental group of the
complement of a link in $S^3$ is obtained from a planar projection). The
generators correspond to the edges of $\Sigma$, indicated by
oriented small undercrossings of the 12 edges of the planar  projection of $\Sigma$.

\medskip

We define a surjection $\psi: \pi_1({\cal O}) \to G$ as follows. If an edge has label
$p$, the image under $\psi$ of the corresponding generator is defined as $\phi(x)$, and
for labels $q, r$ and $s$ the images are respectively $\phi(y), \phi(z)$ and
$\phi(w)$.

\medskip

The orbifold covering of $\cal O$ corresponding to the kernel of $\psi$ is a compact
hyperbolic 3-manifold
$\cal M$ with totally geodesic boundary and an isometric action of $G$. On four of the 
boundary components of $\cal M$ the action of $G$ is equivalent to the given action
on $\cal S$ (on all other boundary components there are actions of cyclic subgroups of
$G$, of orders $p, q, r$ or $s$). By doubling  $\cal M$ and the action of $G$ along the
boundary we obtain a closed hyperbolic 3-manifold with an isometric action of $G$ and a
geodesic embedding of 
$({\cal S},G)$, completing the proof of Theorem 2.

\bigskip  \bigskip

{\bf Remark.}  A proof of Theorem 1 can be given along similar lines, avoiding the
surgery argument in the second part of the proof. Let $({\cal S},G)$ be a 
finite irreducible group action of triangle type $(p,q,r)$, defined by a surjection 
$\phi: (p,q,r) \to  G$. If one of the integers $p, q$ or $r$ is equal to 2, we refer to
the proof given in section 2 using a truncated tetrahedron.

\medskip

If $p,q$ and $r$ are all different from 2, we use a truncated cube instead. 
We consider a planar projection of the 1-skeleton of the cube which consists of two
concentric squares (one contained in the other, with parallel edges) whose vertices
are joined by four further edges (parts of the diagonals of the larger square). Then,
similar as in the proof of Theorem 2, one can label the 12 edges of the cube by $p, q$
and $r$ and choose directions of undercrossings such that, around each vertex, the
labels of the edges are in cyclic order
$p, q, r$. Since we excluded 2 as a label, there are no essential Euclidean
2-suborbifolds (of type (2,2,2,2)). By Andreev's theorem there is a hyperbolic
realization in $\Bbb H^3$ of the orthogonally truncated cube, and the proof is
then similar to the proof of Theorem 2.

\bigskip  \bigskip
\vfill \eject

{\bf 4. Further results}

\medskip

Let $G$ be a finite group of isometries of a hyperbolic surface ${\cal S}$; lifting
$G$ to the universal covering $\Bbb H^2$ of $\cal S$, we get a Fuchsian group $F$ and a
surjection
$\phi: F \to G$ whose kernel is the universal covering group of $\cal S$.

\bigskip

{\bf Theorem 3.}   {\sl For an isometric action $({\cal S},G)$ of a finite group $G$ on 
a hyperbolic surface $\cal S$, suppose that the associated Fuchsian group $F$ is a
subgroup of a triangle group  $(p,q,r)$. Then the isometric action
$({\cal S}, G)$ embeds geodesically into a closed hyperbolic 3-manifold $({\cal M},G)$
with an isometric action of $G$.}

\bigskip

{\it Proof.}  Let ${\cal O}^3$ denote one of the hyperbolic 3-orbifolds constructed in
the proof of Theorem 1 from a truncated hyperbolic tetrahedron (in section 2) or from a
truncated hyperbolic cube  (in the Remark in section 3). Each boundary component of
${\cal O}^3$ is a totally geodesic 2-orbifold ${\cal O}^2$ of triangular type $(p,q,r)$
(a 2-sphere with 3 branch points of orders $p,q$ and $r$), with orbifold fundamental
group the triangle group $(p,q,r)$. Choosing one of these boundary components, it
follows from the proof of Theorem 1 that there is a retraction with torsionfree
kernel of orbifold fundamental groups
$$r: \pi_1({\cal O}^3) \to  \pi_1({\cal O}^2) = (p,q,r)$$  

(a surjective homomorphism
which is the identity on the triangle group $(p,q,r)$). 

\medskip

Let  $\pi_1(\tilde {\cal O}^3) = r^{-1}(F)$ denote the preimage of $F \subset (p,q,r)$
in $\pi_1({\cal O}^3)$, defining  a finite 3-orbifold covering  
$\tilde {\cal O}^3 \to {\cal O}^3$.  Note that the preimage of the totally geodesic
boundary component  ${\cal O}^2$ of  ${\cal O}^3$ is a totally geodesic boundary
component ${\cal F} = \Bbb H^2/F$ of $\tilde {\cal O}^3$ (the covering of ${\cal O}^2$
defined by the subgroup $F \subset (p,q,r) =  \pi_1({\cal O}^2$). 

\medskip

Let $\pi_1(M)$ denote the kernel of  the composition 
$$\psi = \phi \circ r: \pi_1(\tilde {\cal O}^3) \to F \to G,$$
defining a finite regular covering  ${\cal M} \to \tilde {\cal O}^3$, with covering
group $G$ acting by isometries of a compact hyperbolic 3-manifold $\cal M$ with
totally geodesic boundary (note that the kernel of $\psi$ does not contain torsion
elements).  The preimage in $\cal M$  of the totally geodesic boundary component $\cal
F$ of $\tilde {\cal O}^3$  is a totally geodesic boundary component $\cal S$ of $\cal
M$ on which $G$ acts by isometries, realizing the initial action of $G$ on $\cal S$  
associated to the surjection $\phi: F \to G$. 

\medskip

Finally, $({\cal S}, G)$ embeds geodesically into a closed hyperbolic 3-manifold which
is the double of $M$ along the boundary. 

\medskip

This completes the proof of the Theorem 3.

\bigskip
\vfill \eject

Now let  $({\cal S},G)$  be a {\it topological action} (i.e., by homeomorphisms) of  a
finite group $G$ on a surface  $\cal S$ of genus
$g>1$.  As noted in the introduction,
geometrization of $({\cal S},G)$  defines again a Fuchsian group $F$ (unique up to
isomorphism), and Theorem 3 implies:

\bigskip

{\bf Corollary 3.}   {\sl  Let $F$ be a Fuchsian group associated to a topological 
action $({\cal S},G)$ of a finite group $G$ on a closed surface $\cal
S$ of genus $g>1$. If $F$ is isomorphic to a subgroup of a triangle group $(p,q,r)$
then $({\cal S},G)$ embeds geodesically
into a closed hyperbolic 3-manifold, 
for some  geometrization of $({\cal S},G)$.}

\bigskip

Up to isomorphism, a Fuchsian
group $F$ is determined by its signature  $(\bar g; p_1, \ldots, p_n)$ (i.e., the
quotient orbifold $\Bbb H^2/F$ has genus $\bar g$ and $n$ singular points of
degrees $p_1, \ldots, p_n$). This raises naturally the following:

\bigskip

{\bf Question.}  What are the signatures of Fuchsian groups $F$ occuring as 
subgroups of some triangle group (or of some {\it planar Fuchsian group}, considered in
the following)?

\bigskip

Theorem 3 applies to the case of {\it free actions} on surfaces. Let $G$ be a finite
group of fixed-point free homeomorphisms of a closed surface $\cal S$ of genus $g>1$ (a
free topological action). The quotient
${\cal S}/G$ is again a closed surface, and the Fuchsian group $F$  is isomorphic to a
surface group  (the fundamental group of the quotient surface); since all surface
groups are isomorphic to subgroups of triangle groups, Corollary 3 applies. We note
that, in general, a free action of a finite group on a surface does not bound (does not
extend to an action, free or nonfree, of a compact 3-manifold with the
surface as unique boundary component; see [S]);  but, by Corollary 3, at least it embeds
geodesically into a hyperbolic 3-manifold; note that, by the proof of Theorem 3, it
extends to a free action also on the 3-manifold).

\bigskip

{\bf Corollary 4.}   {\sl  Every free topological action $({\cal S},G)$ of a finite
group $G$ on a closed surface $\cal S$ of genus $g>1$ embeds geodesically to a free
action $({\cal M},G)$ of $G$ on a closed hyperbolic 3-manifold $\cal M$.}

\bigskip

Finally we show that topological versions of our main results remain true when we  
the Fuchsian triangle groups are replaced by  {\it planar Fuchsian groups}, i.e.
Fuchsian groups of genus $\bar g = 0$.

\bigskip

{\bf Theorem 4.}   {\sl  let $({\cal S}, G)$ be a topological action of a finite group
$G$ on a closed surface
$\cal S$ of genus $g \ge 2$. Suppose that the Fuchsian group $F$ associated to the
action (unique up to isomorphism)  is planar, with signature $(0; p_1, \ldots, p_n)$
different from $(0;2,2, \ldots, 2)$. Then, 
for some geometrization of $({\cal S}, G)$, the action embeds geodesically into a
closed hyperbolic 3-manifold 
$({\cal M}, G)$.}

\bigskip

{\it Proof.}  Suppose $p_1 > 2$. Let $\phi:F \to G$ be the surjection associated to the
action.  The proof is along similar lines as the proofs of
Theorems 1 and 2.

\medskip

Let  ${\cal O}^3$ be a hyperbolic 3-orbifold  with totally geodesic boundary whose
singular set consists of the 1-skeleton of an orthogonally truncated  polyhedron which
is the  double cone over an $n$-gon (an octahedron if $n =  4$).  Considering a planar
projection of its 1-skeleton (with one cone point at infinity), at each of the two
truncated cone points we associate labels $p_1, \ldots, p_n$ in cyclic order to the
$n$ incoming edges, and labels $p_1, p_i, p_1, p_i$ to the four edges at the $i'th$
vertex on the central $n$-gon
(since   $p_1 > 2$, this avoids Euclidean vertices of type (2,2,2,2)).  Let
${\cal O}^2$ be the 2-sphere with $n$ singular points of orders
$p_1, \ldots, p_n$; note that this is the 2-orbifold occuring as the totally geodesic
boundary component of
${\cal O}^3$ at each of the two cone points. Choosing coherent orientations of
the generators of $\pi_1({\cal O}^3)$ (represented by directed undercrossings of the
$3n$ edges), it is easy to see that there is a retraction of orbifold fundamental
groups $r: \pi_1({\cal O}^3) \to  \pi_1({\cal O}^2)$, and then the  proof finishes 
as the proof of Theorem 3.

\bigskip

Similar as in the proof of Theorem 3, this implies:

\bigskip

{\bf Corollary 5.}   {\sl  If the Fuchsian group  $F$ associated to a topological 
action  $({\cal S}, G)$  is isomorphic to a subgroup of a planar
Fuchsian group not of type $(0;2,2, \ldots, 2)$ then, for some geometrization
of $({\cal S}, G)$, the action embeds geodesically into a closed  hyperbolic
3-manifold.}

\bigskip

\bigskip \bigskip

\centerline {\bf References}

\bigskip

\item {[A]}  E.M. Andreev,  {\it On convex polyhedra in Lobachevsky space.} 
Math. USSR Sb. 10  (1970),  413-440

\item {[B]}  L.A. Best,  {\it  On torsion-free discrete subgroups of PSL(2, \Bbb C)
with compact orbit space.} Can. J. Math.  23  (1971), 451-460

\item {[GZ]}  M. Gradolato, B. Zimmermann,  {\it  Extending finite group actions on
surfaces to hyperbolic 3-manifolds.}  Math. Proc. Camb. Phil. Soc. 117  (1995),
137-151

\item {[S]}  E.G. Samperton,  {\it  Free actions on surfaces that do not extend to
arbitrary actions on 3-manifolds. }  Comptes Rendus Math. 360 (2022), 161-167

\item {[SW]}  J.-C. Schlage-Puchta, J. Wolfart,  {\it  How many quasiplatonic
surfaces?} Archiv Math.  86  (2006), 129-132

\item {[T]} W. Thurston,  {\it The geometry and topology of 3-manifolds.} 
Lecture Notes, Dept. of Math. Princeton Univ. 1977

\item {[V]} E.B. Vinberg,  {\it Geometry II.}  Encyclopaedia Math. Sciences vol. 29,
Springer 1993

\item {[WZ]} C. Wang, S. Wang, Y. Zhang, B. Zimmermann,  {\it Finite group
actions on the genus-2 surface, geometric generators and extendability.}  
Rend. Istit. Mat. Univ. Trieste 52 (2020), 513-524

\item {[Z1]} B. Zimmermann,  {\it Hurwitz groups, maximal reducible groups and maximal
handlebody groups.}  arXiv:2110.11050

\item {[Z2]}  B. Zimmermann, {\it  A note on large bounding and non-bounding
finite group-actions on surfaces of small genus.}  arXiv:2205.14425

\item {[Z3]}  B. Zimmermann, {\it  A note on surfaces bouonding hyperbolic
3-manifolds.}  Monatsh. Math. 142  (2004), 267-273

\item {[Z4]}  B. Zimmermann, {\it  Hurwitz groups and finite group actions on
hyperbolic 3-manifolds.}  J. London Math. Soc. 52  (1995), 199-208

\bye